\documentclass[11pt]{amsart}

\usepackage[utf8]{inputenc}
\usepackage[T1]{fontenc}
\usepackage{amsmath}
\usepackage{amssymb}
\usepackage{mathtools}
\usepackage{enumitem}
\usepackage{xcolor}
\usepackage[pagebackref]{hyperref}
\usepackage[capitalize]{cleveref} 

\hypersetup{
	linktoc=page,
	colorlinks=true,
	linkcolor={red!90!black},
	citecolor={blue!90!black},
	urlcolor=magenta,
	filecolor=teal,
}

\renewcommand*{\backref}[1]{}
\renewcommand*{\backrefalt}[4]{%
	\ifcase #1 (Not cited.)%
	\or        (Cited on page~#2)%
	\else      (Cited on pages~#2)%
	\fi}

\crefname{equation}{}{} 
\crefname{page}{p.}{pp.}
\crefname{enumi}{item}{items}

\newtheorem{theorem}{Theorem}[section]
\newtheorem{lemma}[theorem]{Lemma}

\newtheorem*{claim*}{Claim}

\newtheorem*{theorem*}{Theorem}

\theoremstyle{definition}

\crefname{theorem}{Theorem}{Theorems}
\crefname{lemma}{Lemma}{Lemmas}
\crefname{mlemma}{Main Lemma}{Main Lemmas}
\crefname{keylemma}{Key Lemma}{Key Lemmas}
\crefname{coro}{Corollary}{Corollaries}
\crefname{propo}{Proposition}{Propositions}
\crefname{mpropo}{Main Proposition}{Main Propositions}
\crefname{claim}{Claim}{Claims}
\crefname{definition}{Definition}{Definitions}
\crefname{example}{Example}{Examples}
\crefname{xca}{Exercise}{Exercises}
\crefname{rem}{Remark}{Remarks}

\AddToHook{env/lemma/begin}{\crefalias{theorem}{lemma}}
\AddToHook{env/mlemma/begin}{\crefalias{theorem}{mlemma}}
\AddToHook{env/keylemma/begin}{\crefalias{theorem}{keylemma}}
\AddToHook{env/coro/begin}{\crefalias{theorem}{coro}}
\AddToHook{env/propo/begin}{\crefalias{theorem}{propo}}
\AddToHook{env/mpropo/begin}{\crefalias{theorem}{mpropo}}
\AddToHook{env/claim/begin}{\crefalias{theorem}{claim}}
\AddToHook{env/definition/begin}{\crefalias{theorem}{definition}}
\AddToHook{env/example/begin}{\crefalias{theorem}{example}}
\AddToHook{env/xca/begin}{\crefalias{theorem}{xca}}
\AddToHook{env/rem/begin}{\crefalias{theorem}{rem}}

\numberwithin{equation}{section}

\newcommand{\R}{{\mathbb R}}
\newcommand{\Z}{{\mathbb Z}}
\newcommand{\BZ}{{\mathcal B}}
\newcommand{\EE}{{\mathcal E}}
\newcommand{\HH}{{\mathcal H}}
\newcommand{\M}{{\mathcal M}}
\newcommand{\PP}{{\mathcal P}}
\newcommand{\characteristic}{\mathbf{1}} 
\newcommand{\Capacity}{\operatorname{Cap}}
\newcommand{\divv}{\operatorname{div}}
\newcommand{\loc}{\operatorname{loc}}

\begin{document}

\title{On the packing-dimension analogue of Bourgain's constant}

\author{Ignasi Guill{\'e}n-Mola}
\address{Ignasi Guill{\'e}n-Mola}
\email{ignasiguillen@gmail.com}

\date{\today}



\begin{abstract}
    For every $n\geq 2$, we construct a domain in $\mathbb{R}^n$ whose harmonic measure has packing dimension $n$. That is, the packing analogue of Bourgain's constant is zero, in contrast with Bourgain's constant, which is positive. Moreover, the elliptic measure of the same domain also has packing dimension $n$ for every uniformly elliptic matrix.
\end{abstract}

\maketitle


\section{Introduction}

It is a classical, long-standing open problem in harmonic analysis and geometric measure theory to determine the precise gap between the Hausdorff dimension of the harmonic measure $\omega_\Omega$ of a domain (open and connected) $\Omega\subset \R^n$ and the dimension of the ambient space $\R^n$. That is, to find the value of the so-called Bourgain constant
$$
b_n \coloneqq n - \sup_{\Omega\subset \R^n} \dim_\HH \omega_\Omega.
$$
Below we recall the definition of harmonic measure, and here we just note that the Hausdorff dimension of $\omega_\Omega$ is
$$
\dim_\HH \omega_\Omega = \inf \{\dim_\HH E : \omega_\Omega(E^c)=0\}.
$$

In this note we study the analogue of Bourgain's constant for the packing dimension. Given a set $E\subset\R^n$ and $s\geq 0$, the $s$-dimensional packing premeasure of $E$ is
$$
\PP^s_0(E) \coloneqq \lim_{\delta\to 0} \sup\Big\{\sum_{i\geq 1} (2r_i)^s : \{B(x_i,r_i)\}_{i\geq 1} \text{ disjoint},\ x_i\in E,\ r_i\leq \delta\Big\},
$$
the $s$-dimensional packing measure of $E$ is
$$
\PP^s(E) \coloneqq \inf\Big\{\sum_{i\geq 1} \PP^s_0(E_i) : E\subset \bigcup_{i\geq 1} E_i\Big\},
$$
and the packing dimension of $E$ is defined as
$$
\dim_P E \coloneqq \inf\{s\geq 0 : \PP^s(E)=0\} = \sup\{s\geq 0: \PP^s(E)=\infty\}.
$$
As for the Hausdorff dimension, the packing dimension of $\omega_\Omega$ is
$$
\dim_P \omega_\Omega = \inf \{\dim_P E : \omega_\Omega (E^c)=0\},
$$
and the packing analogue of Bourgain's constant is
$$
p_n \coloneqq n - \sup_{\Omega\subset\R^n} \dim_P \omega_\Omega.
$$
We study the behavior of $p_n$ even for elliptic measure, which we define below.

Let $\Omega\subset \R^n$, $n\geq 2$, be a bounded domain and let $A\in L^\infty (\Omega; \R^{n\times n})$ be a (not necessarily symmetric) uniformly elliptic matrix. That is, a real and measurable matrix $A:\Omega\to \R^{n\times n}$ such that there exists $\lambda\geq 1$, the so-called ellipticity constant, for which
$$
\lambda^{-1}|\xi|^2\leq \langle A(x)\xi, \xi \rangle \text{ and } \langle A(x)\xi,\eta\rangle\leq \lambda|\xi||\eta| \text{ for all } \xi,\eta\in\R^n \text{ and a.e.\ } x\in\Omega.
$$
We consider the second-order operator in divergence form $L_A \coloneqq -\divv (A\nabla \cdot)$, understood in the weak sense. We say that a function $u\in W^{1,2}_{\loc}(\Omega)$ is $L_A$-harmonic in $\Omega$ if $L_A u=0$ in $\Omega$, that is,
$$
\int_\Omega \langle A\nabla u, \nabla \varphi\rangle \, dx = 0 \text{ for all } \varphi \in C^\infty_c(\Omega).
$$
By the De Giorgi--Nash--Moser theory, $L_A$-harmonic functions are locally H\"older continuous and satisfy Harnack's inequality and the maximum principle.

Perron's method provides, for every $f\in C(\partial\Omega)$, an $L_A$-harmonic function $h_f^A$ in $\Omega$, the so-called Perron solution with boundary data $f$. By the maximum principle and the Riesz representation theorem, for every $x\in\Omega$ there exists a Radon probability measure $\omega^x_{\Omega,A}$ on $\partial\Omega$ such that
$$
h_f^A(x) = \int_{\partial\Omega} f\,d\omega^x_{\Omega,A} \text{ for all } f\in C(\partial\Omega).
$$
The measure $\omega^x_{\Omega,A}$ is the elliptic measure of $\Omega$ with pole at $x$, and when $A$ is the identity matrix it is the harmonic measure $\omega^x_\Omega$ of $\Omega$ with pole at $x$. For a detailed construction of elliptic measure and further properties, see \cite[Chapter 11]{Heinonen2006}. By Harnack's inequality, the measures $\omega^x_{\Omega,A}$ and $\omega^y_{\Omega,A}$ are mutually absolutely continuous for all $x,y\in\Omega$, so the Hausdorff and packing dimensions of $\omega^x_{\Omega,A}$ do not depend on the pole $x$. For this reason, we drop the pole from the notation and write $\omega_{\Omega,A}$ (and $\omega_\Omega$ when $A$ is the identity) whenever the pole plays no role.

In the planar case $n=2$, {\O}ksendal \cite{Oksendal1981} showed that $\omega_\Omega$ is singular with respect to area and conjectured that $\dim_\HH\omega_\Omega\leq 1$ for every domain $\Omega\subset \R^2$. The conjecture was first verified for broad classes of domains: Carleson \cite{Carleson1985} proved that $\dim_\HH\omega_\Omega\leq 1$ for ``snowflake type'' domains and $\dim_\HH\omega_\Omega<1$ for the complement of some self-similar ``Cantor type'' sets, and Makarov \cite{Makarov1985} showed that $\dim_\HH\omega_\Omega=1$ for every simply connected domain. Jones and Wolff \cite{Jones1988} proved {\O}ksendal's conjecture for all planar domains, that is, $b_2 \geq 1$, which turns out to be optimal by taking the unit ball, say. Years later, Wolff \cite{Wolff1993} refined this result by showing that the harmonic measure in the plane is concentrated on a set of $\sigma$-finite $\HH^1$-measure.

In higher dimensions $n\geq 3$, where {\O}ksendal's singularity result also holds, Bourgain \cite{Bourgain1987} improved it by proving that $b_n>0$, and Wolff \cite{Wolff1995} constructed a domain in $\R^3$ whose harmonic measure has Hausdorff dimension greater than $2$ (see also \cite{Lewis-Verchota-Vogel-2005} for higher dimensions), whence $b_n < 1$. As mentioned above, the optimal value $b_n$ is still unknown for $n\geq 3$, whereas $b_2 = 1$, and Bishop \cite{Bishop1992} conjectured that $b_n = 1/(n-1)$. Badger and Genschaw \cite{Badger2024} refined Bourgain's proof to find the first quantitative lower bounds for $b_n$, namely $b_n\geq c\, n^{-2n(n-1)}/\log n$ for all $n\geq 3$, where $c>0$ is a constant independent of $n$, and $b_3\geq 10^{-15}$, that is, $\dim_\HH \omega_\Omega \leq 2.99999\,99999\,99999$ for every domain $\Omega\subset \R^3$. Very recently, Albert, Badger and Genschaw \cite{Albert-Badger-Genschaw-2026} improved both bounds, obtaining $b_n\geq (n-1)^{-2n-6.12\sqrt{n}}$ for all $n\geq 2$, so that the order of magnitude goes from $n^{-O(n^2)}$ to $n^{-O(n)}$, and $b_3> 10^{-5}$, that is, $\dim_\HH \omega_\Omega < 2.99999$ for every domain $\Omega\subset \R^3$. These bounds are still far from the conjecture, if true.

In \cite{Albert-Badger-Genschaw-2026} the authors also propose ten challenges\footnote{The author learned of these challenges some time before \cite{Albert-Badger-Genschaw-2026} appeared.} to better understand the dimension of harmonic measure. Challenge 9, on the dimension drop for harmonic measure on the complement of Ahlfors regular sets in the plane, has been very recently solved by Cai, Jeznach and Tolsa \cite{Cai-Jeznach-Tolsa-close-to-1,Cai-Jeznach-Tolsa-flat}. Challenge 10 concerns elliptic measure, which we also consider here. In this case, Sweezy showed that the Hausdorff dimension of elliptic measure can be arbitrarily close to $n$ if the ellipticity constant is large enough, both in the plane \cite{Sweezy1992} and in higher dimensions \cite{Sweezy1994}. In the plane, the author \cite{Guillen-quasiconformal} proved that $\dim_\HH\omega_{\Omega,A}\leq 2K_\lambda/(K_\lambda+1)$, where $K_\lambda=\lambda$ if $A$ is symmetric and $K_\lambda=\lambda^2+\sqrt{\lambda^4-1}$ otherwise, and that $\dim_\HH\omega_{\Omega,A}\leq 1$ if $A$ is continuous up to the boundary. Moreover, the author, Prats and Tolsa \cite{Guillen-Prats-Tolsa} showed that $\omega_{\Omega,A}$ is concentrated on a set of $\sigma$-finite $\HH^1$-measure if $\Omega$ is a planar Reifenberg flat domain with small constant and $A$ is Lipschitz. In higher dimensions $n\geq 3$, it is expected that Bourgain's theorem extends to elliptic measure, that is, that $\dim_\HH\omega_{\Omega,A}\leq n-b$ for some $b>0$ depending only on $n$ and $\lambda$, as already remarked by {\O}ksendal \cite[p.~295]{Oksendal1988} and in \cite{Albert-Badger-Genschaw-2026}. Indeed, the key ingredient in Bourgain's proof is a lower bound for harmonic measure in terms of Hausdorff content \cite[Lemma 1]{Bourgain1987}, known as Bourgain's alternative, whose proof also works for $L_A$, with constants depending on $\lambda$. Thus, the main difficulty in Challenge 10 is to understand how the optimal value $b_n(\lambda)$ of such $b$ depends on the ellipticity constant $\lambda$.

Challenge 7, which we answer in this note, asks whether $p_n>0$ for every $n\geq 2$ and, if so, whether $p_n=b_n$. By the definitions of Hausdorff and packing dimensions, there holds $p_n\leq b_n$, as we always have $\dim_\HH E\leq \dim_P E\leq n$ for all $E\subset \R^n$. As pointed out in \cite{Albert-Badger-Genschaw-2026}, the evidence so far was mixed. On the one hand, Housworth \cite{Housworth1995} proved the packing analogue of Makarov's result: the harmonic measure of every simply connected planar domain is concentrated on a set of packing dimension $1$. On the other hand, the lower density conjecture, proved by Choi \cite{Choi2004}, implies (by \cite[Theorem 6.11]{Mattila1995}) that the harmonic measure of many planar domains, such as the von Koch snowflake domain, vanishes on every set of $\sigma$-finite $\PP^1$-measure, unlike in the Hausdorff case \cite{Wolff1993}.

Moreover, we note that by {\O}ksendal's theorem (which also holds for elliptic measure \cite[Theorem 5.9]{Oksendal1988}) and $\PP^n(F)\approx \HH^n(F)$ for all $F\subset\R^n$ (see \cite[Theorem 6.11]{Mattila1995}), $\omega_{\Omega,A}$ is singular with respect to $\PP^n$ for every bounded domain $\Omega\subset\R^n$. That is, there exists a set $E\subset\partial\Omega$ such that
\begin{equation}\label{eq:packing null}
\omega_{\Omega,A}(E^c)=0 \text{ and } \PP^n(E)=0.
\end{equation}
This might suggest that $\dim_P\omega_{\Omega,A}<n$ for every domain $\Omega$, as is the case for the Hausdorff dimension. In contrast, here we show that the packing analogue of Bourgain's constant is zero, even for elliptic measure.

\begin{theorem}\label{thm:main}
    There exists a bounded domain $\Omega\subset\R^n$ such that $\dim_P \omega_{\Omega, A} = n$ for all uniformly elliptic matrices $A$.
\end{theorem}

In particular, taking $A$ to be the identity matrix, we get $p_n=0$ for all $n\geq 2$, which answers Challenge 7 in \cite{Albert-Badger-Genschaw-2026} in the negative. Note that the same domain $\Omega$ works for all uniformly elliptic matrices $A$, as the ellipticity constant only affects the scale from which our estimates hold.

By \eqref{eq:packing null}, \cref{thm:main} is sharp at the level of measures: although $\dim_P\omega_\Omega$ can be $n$, $\omega_\Omega$ is always concentrated on a set of zero $\PP^n$-measure. This is in contrast with the planar case, where $\omega_\Omega$ is always concentrated on a set of $\sigma$-finite $\HH^1$-measure by Wolff's theorem \cite{Wolff1993}, and in general not on a set of finite $\HH^1$-measure. For packing measures, the lower density conjecture implies that $\omega_\Omega$ need not be concentrated on a set of $\sigma$-finite $\PP^1$-measure (see above), and \cref{thm:main} shows that the same holds for $\PP^s$ for every $s<2$, whereas $\omega_\Omega$ is always concentrated on a set of zero $\PP^2$-measure.

The proof relies on the difference between the Hausdorff and the packing dimensions of a measure $\mu$. Roughly speaking, the former is determined by the lower local dimension of $\mu$, and the latter by the upper local dimension of $\mu$, that is, by
$$
\liminf_{r\to 0}\frac{\log\mu(B(x,r))}{\log r} \text{ and } \limsup_{r\to 0}\frac{\log\mu(B(x,r))}{\log r},
$$
respectively (see \cite[Proposition 2.3]{Falconer1997}). Hence, to bound the Hausdorff dimension of $\mu$ from below one needs an estimate of the form $\mu(B(x,r))\lesssim r^s$ at all small scales $r$, whereas for the packing dimension it suffices to have it along a sequence of scales tending to zero, which gives much more room. For instance, by Bourgain's theorem, for every $s>n-b_n$ and $\omega_\Omega$-almost every $x$ there are arbitrarily small radii $r$ with $\omega_\Omega(B(x,r))\geq r^s$, but this does not prevent $\omega_\Omega(B(x,r))$ from being much smaller at other scales.

Our domain is $\Omega=B(0,1)\setminus K$, where the compact set $K$ exploits this freedom. It is a Cantor-type set with scales $r_k\ll\delta_k\ll r_{k-1}$ that become increasingly lacunary as $k$ increases (see \cref{sec:construction}). At stage $k\geq 1$, each ball of the previous generation, which has radius $r_{k-1}$, is replaced by many tiny balls of radius $r_k$ centered at the points of a $\delta_k$-net of it. Since $r_k\ll\delta_k$, every ball of radius $\delta_k/2$ centered at $K$ intersects $K$ only inside one of these tiny balls, so its elliptic measure is at most that of a ball of radius $r_k$. Exploiting the huge gap between consecutive scales, we show that the latter is at most a constant times $\delta_k^{n-o(1)}$. This forces the upper local dimension of the elliptic measure to be (at least) $n$ at every point of $K$. Moreover, the radii $r_k$ are chosen so that $K$ has positive capacity, which guarantees that $K$ has positive elliptic measure.

The paper is organized as follows. In \cref{sec:construction} we construct the compact set $K$. In \cref{sec:elliptic measure} we show that the elliptic measure of every ball of radius $\delta_k/2$ centered at $K$ is at most a constant times $\delta_k^{n-2/(k+1)}$ (\cref{lemma:upper bound}) and that $K$ has positive elliptic measure (\cref{lemma:positive elliptic measure}). Finally, in \cref{sec:proof} we prove \cref{thm:main}.

\subsection*{Notation}

Throughout the paper $n\geq 2$. We write $a\lesssim b$ if there exists a constant $C\geq 1$ such that $a\leq Cb$, $a\approx b$ if $a\lesssim b\lesssim a$, and $a\ll b$ if $a\leq c\, b$ for a sufficiently small constant $c>0$. All the constants are allowed to depend on $n$ and on the ellipticity constant $\lambda$ unless explicitly stated otherwise: if they depend only on a parameter $\kappa$, we write $a\lesssim_\kappa b$, $a\approx_\kappa b$ and $a\ll_\kappa b$. Balls are open, $B(x,r) \coloneqq \{y\in\R^n:|y-x|<r\}$. Given a set $E$, we write $\# E$ for its cardinality, $\overline{E}$ for its closure and $E^c\coloneqq\R^n\setminus E$. Once the domain $\Omega$ is fixed in \cref{sec:construction}, we write $\omega^p_A\coloneqq\omega^p_{\Omega,A}$ and $\omega_A\coloneqq\omega_{\Omega,A}$.

\section{Construction of \texorpdfstring{$K$}{K}}\label{sec:construction}

Fix $0<r_0\ll_n 1/4$. For $k\geq 1$, let
$$
\delta_k \coloneqq r_{k-1}^{k+1},\quad
N_k \coloneqq \# \left(\delta_k \Z^n \cap B(0,r_{k-1}/2)\right),\quad
r_k \text{ such that }
\begin{cases}
    N_k r_k^{n-2} = 2 r_{k-1}^{n-2} & \text{if }n\geq 3,\\
    \frac{N_k}{\log(1/r_k)} = \frac{2}{\log(1/r_{k-1})} &\text{if }n=2.
\end{cases}
$$
It is clear that $r_k\leq r_{k-1}/2$ for $k\geq 1$, since $N_k$ is large because $r_0\ll_n 1/4$. In particular, the sequence $\{r_k\}_{k\geq 0}$ is strictly decreasing and converges to zero.

Let $\BZ_0 \coloneqq \{B(0,r_0)\}$, and inductively,
$$
\BZ_k \coloneqq \{ B(z,r_k) : B(c,r_{k-1})\in \BZ_{k-1} \text{ and } z\in (c+\delta_k\Z^n)\cap B(c,r_{k-1}/2)\}.
$$
From now on we fix the domain
$$
\Omega \coloneqq B(0,1) \setminus K, \text{ where } K\coloneqq \bigcap_{k\geq 1} \bigcup_{B\in \BZ_k} \overline{B}.
$$
Note that $K$ is compact, since it is an intersection of finite unions of closed balls. Moreover, each ball of $\BZ_k$ is contained in the ball of $\BZ_{k-1}$ from which it is constructed, since $r_k\leq r_{k-1}/2$, and every center of a ball of $\BZ_k$, $k\geq 0$, belongs to $K$, since $0\in\delta_j\Z^n$ for all $j$.

Now we claim that
\begin{equation}\label{eq:rk small}
r_k \leq \delta_k/4,
\end{equation}
which in particular implies that $\Omega$ is connected. Indeed, using
$$
N_k \approx \left(\frac{r_{k-1}}{\delta_k}\right)^n,
$$
from the definition of $N_k$, if $n\geq 3$ we have
$$
\left(\frac{r_k}{\delta_k}\right)^{n-2} = \frac{2 r_{k-1}^{n-2}}{N_k \delta_k^{n-2}} \leq C \frac{\delta_k^2}{r_{k-1}^2} = C r_{k-1}^{2k} \leq C r_0^2 \leq (1/4)^{n-2}
$$
as claimed, where we used in the last step that $r_0\ll_n 1/4$; and if $n=2$, since $\log(1/r_{k-1})\geq 1$ and $r_{k-1}=\delta_k^{1/(k+1)}$, we have
$$
\frac{1}{\log(1/r_k)} = \frac{2}{N_k\log(1/r_{k-1})} \leq C \frac{\delta_k^2}{r_{k-1}^2} = C \delta_k^{\frac{2k}{k+1}} \leq C\delta_k \leq \frac{1}{\log(4/\delta_k)},
$$
where we used in the last step that $\delta_k\leq r_0^2$ and $r_0\ll_n 1/4$. That is, $r_k\leq \delta_k/4$, as claimed.

\section{Elliptic measure of \texorpdfstring{$K$}{K}}\label{sec:elliptic measure}

\subsection{Upper bound at every scale}

In this subsection we prove the following upper bound for points in $K$. This is essentially what we need to prove \cref{thm:main}, modulo showing that $K$ is not invisible for the elliptic measure, that is, $\omega_A(K)>0$, which we do in \cref{lemma:positive elliptic measure} below.

\begin{lemma}\label{lemma:upper bound}
    For all $k\geq 1$, $x\in K$ and $p\in\Omega$ with $|p|=1/2$, there holds $\omega^p_A (B(x,\delta_k / 2))\lesssim \delta_k^{n-2/(k+1)}$.
    \begin{proof}
        Extend $A$ by the identity outside $\Omega$, and let $B(z,r_k)\in \BZ_k$ be the ball containing $x$. By \eqref{eq:rk small}, we have that $\omega_A^p (B(x,\delta_k / 2)) = \omega_A^p (B(z,r_k))$. By the maximum principle, as $\Omega\subset B(z,2)$,\footnote{When $A$ is the identity, the last term below is explicit: it equals $(|p-z|^{2-n}-2^{2-n})/(r_k^{2-n}-2^{2-n})$ if $n\geq 3$ and $\log(2/|p-z|)/\log(2/r_k)$ if $n=2$, which already gives the estimate below, since $1/4<|p-z|<3/4$.}
        $$
        \omega_A^p (B(z,r_k))
        \leq \omega^p_{\Omega\setminus \overline{B(z,r_k)},A} (\overline{B(z,r_k)})
        \leq \omega^p_{B(z,2)\setminus \overline{B(z,r_k)},A} (\overline{B(z,r_k)}).
        $$
        Let $\EE^A_z$ be the elliptic fundamental solution of $L_A$ with pole at $z$, constructed by Kenig and Ni in \cite[Appendix]{Kenig1985} in the plane and by Hofmann and Kim in \cite[Theorem 3.1]{Hofmann2007} in higher dimensions. Since $\EE^A_z$ is $L_A$-harmonic in $\R^n\setminus\{z\}$ and $\EE^A_z(\xi)\to\infty$ as $\xi\to z$, for $k$ large enough the function
        $$
        u(\xi)\coloneqq\frac{\EE^A_z(\xi)-\min_{\partial B(z,2)}\EE^A_z}{\min_{\partial B(z,r_k)}\EE^A_z - \min_{\partial B(z,2)}\EE^A_z}
        $$
        is well defined, $L_A$-harmonic in $B(z,2)\setminus\overline{B(z,r_k)}$, nonnegative on $\partial B(z,2)$ and at least $1$ on $\partial B(z,r_k)$. Hence, by the maximum principle,
        $$
        \omega^p_{B(z,2)\setminus \overline{B(z,r_k)},A} (\overline{B(z,r_k)}) \leq u(p).
        $$
        As $\EE^A_z(\xi) \gtrsim |z-\xi|^{2-n}$ if $n\geq 3$ and $\EE^A_z(\xi) \gtrsim \log (1/|z-\xi|)$ if $n=2$ as $\xi\to z$ (by comparison with the Green function \cite[(1.9)]{Gruter1982} and by \cite[Theorem A-2]{Kenig1985}, respectively), and the numerator of $u(p)$ is bounded by a constant since $|p|=1/2$ (by \cite[(3.55)]{Hofmann2007} and \cite[Theorem A-2]{Kenig1985}, respectively), we get that there is $k_0=k_0(\lambda)$ such that if $k\geq k_0$, then
        $$
        \omega^p_{B(z,2)\setminus \overline{B(z,r_k)},A} (\overline{B(z,r_k)})\lesssim
        \begin{cases}
            r_k^{n-2}, &\text{if }n\geq 3,\\
            \frac{1}{\log(1/r_k)}, &\text{if }n=2.
        \end{cases}
        $$
        If $n\geq 3$, by $N_k r_k^{n-2} = 2r_{k-1}^{n-2}$ and $\delta_k = r_{k-1}^{k+1}$ in the definition of $K$, and $N_k \approx (r_{k-1}/\delta_k)^n = \delta_k^{-nk/(k+1)}$, we have
        $$
        r_k^{n-2} = \frac{2r_{k-1}^{n-2}}{N_k} \approx \delta_k^{\frac{n-2}{k+1}}\,\delta_k^{\frac{nk}{k+1}} = \delta_k^{n-\frac{2}{k+1}}.
        $$
        Similarly, if $n=2$, since $\log(1/r_{k-1})\geq 1$,
        $$
        \frac{1}{\log(1/r_k)} = \frac{2}{N_k\log(1/r_{k-1})} \leq \frac{2}{N_k} \approx \delta_k^{\frac{2k}{k+1}} = \delta_k^{2-\frac{2}{k+1}}.
        $$
        This proves the lemma for $k\geq k_0$. For $k<k_0$ it is trivial, since $\|\omega_A^p\|\leq 1$ and $\delta_k\geq\delta_{k_0}$, which depends only on $n$ and $\lambda$.
    \end{proof}
\end{lemma}

\subsection{Positivity of the elliptic measure of \texorpdfstring{$K$}{K}}

In this subsection we aim to prove the following lemma.

\begin{lemma}\label{lemma:positive elliptic measure}
    The set $K$ satisfies $\omega_A (K)>0$.
\end{lemma}

Given a bounded domain $D$, a point $x\in \partial D$ is called Wiener regular for $L_A$ if $h_f^A(y)\to f(x)$ as $D\ni y\to x$ for every $f\in C(\partial D)$, where $h_f^A$ is the Perron solution in $D$ with boundary data $f$. By a result of Littman, Stampacchia and Weinberger \cite{Littman-Stampacchia-Weinberger-1963}, the set of Wiener regular points of $D$ does not depend on the operator $L_A$, so we will use the theory for the harmonic case. This theory is formulated in terms of capacity, which we now define.

First we define, for $t>0$,
$$
\EE (t) =
\begin{cases}
   \dfrac{t^{2-n}}{(n-2)\,\HH^{n-1}(\mathbb{S}^{n-1})}, &\text{if }n\geq 3,\\
   \dfrac{1}{2\pi}\log\dfrac{1}{t}, &\text{if }n=2,
\end{cases}
$$
so that $\EE(|x|)$ is the fundamental solution of $-\Delta$ in $\R^n$ with pole at $0$. Given a bounded set $E$ and a Radon probability measure $\mu$ on $E$ (we write $\mu \in \M_1 (E)$), we define the potential
$$
U_\mu (x) = \int_E \EE(|x-y|) \,d \mu (y) \text{ for }x \in E,
$$
and the energy
$$
I(\mu) = \int_E U_\mu(x) \, d\mu(x).
$$
The capacity is defined as
$$
\Capacity (E)= \frac{1}{\inf_{\mu \in \M_1(E)} I(\mu)}.
$$

In order to prove \cref{lemma:positive elliptic measure}, we first need to see that the set $K$ has positive capacity.

\begin{lemma}
The set $K$ satisfies $\Capacity (K)>0$.
\begin{proof}
By definition, it suffices to get a uniform upper bound for $U_\mu (x)$ for all $x\in K$ for a measure $\mu \in \M_1(K)$.

Let $\mu$ be the probability measure supported on $K$ such that $\mu(K\cap\overline{B}) = m_k\coloneqq (N_1\cdots N_k)^{-1}$ for all $B\in\BZ_k$ and $k\geq 1$, and let $m_0\coloneqq 1$. Let $x\in K$. For $j\geq 0$, let $B^j(x)\in\BZ_j$ be the ball whose closure contains $x$, and for $j\geq 1$ let $K^j\coloneqq (K\cap\overline{B^{j-1}(x)})\setminus\overline{B^j(x)}$. Since $\mu(\{x\})\leq m_j\to 0$ as $j\to\infty$, we have
$$
U_\mu(x) = \sum_{j\geq 1}\int_{K^j}\EE(|x-y|)\,d\mu(y).
$$

Fix $j\geq 1$, and let $B^{j-1}(x)=B(c,r_{j-1})$ and $B^j(x)=B(z',r_j)$. Then $K^j$ is contained in the union of the closures of the balls $B(z'',r_j)\in\BZ_j$ with $z''\in(c+\delta_j\Z^n)\cap B(c,r_{j-1}/2)$ and $z''\neq z'$. For such $z''$ and $y\in\overline{B(z'',r_j)}$, since $z'-z''\in\delta_j\Z^n\setminus\{0\}$ and $r_j\leq\delta_j/4$ by \eqref{eq:rk small}, we have
$$
|x-y|\geq |z'-z''|-2r_j\geq\frac{|z'-z''|}{2},
$$
and $|z'-z''|<r_{j-1}$ because $z',z''\in B(c,r_{j-1}/2)$. Since $\EE$ is decreasing and each of these balls has measure $m_j$, we get
$$
\int_{K^j}\EE(|x-y|)\,d\mu(y)\leq m_j\sum_{\substack{v\in\delta_j\Z^n\setminus\{0\}\\ |v|<r_{j-1}}}\EE(|v|/2).
$$
To estimate the sum, note that the balls $B(v,\delta_j/2)$, $v\in\delta_j\Z^n$, are pairwise disjoint and, if $|v|<r_{j-1}$, contained in $B(0,2r_{j-1})$. Moreover, for $v\neq 0$ and $w\in B(v,\delta_j/2)$ we have $|w|\leq|v|+\delta_j/2\leq 3|v|/2$, so $\EE(|v|/2)\leq\EE(|w|/3)$. Hence, integrating in polar coordinates and using that $r_{j-1}\leq r_0\ll_n 1/4$,
$$
\sum_{\substack{v\in\delta_j\Z^n\setminus\{0\}\\ |v|<r_{j-1}}}\EE(|v|/2)\lesssim\frac{1}{\delta_j^n}\int_{B(0,2r_{j-1})}\EE(|w|/3)\,dw\lesssim\left(\frac{r_{j-1}}{\delta_j}\right)^n\EE(r_{j-1})\approx N_j\,\EE(r_{j-1}).
$$
Since $m_jN_j=m_{j-1}$, we obtain
$$
\int_{K^j}\EE(|x-y|)\,d\mu(y)\lesssim m_{j-1}\EE(r_{j-1}).
$$

Finally, the constants in the definition of $\EE$ cancel out in the definition of $r_k$, which reads $N_k/\EE(r_k)=2/\EE(r_{k-1})$ both if $n\geq 3$ and if $n=2$. Thus $m_k\EE(r_k)=m_{k-1}\EE(r_{k-1})/2$ for all $k\geq 1$, that is, $m_k\EE(r_k)=2^{-k}\EE(r_0)$, and so
$$
U_\mu(x)\lesssim\sum_{j\geq 1}2^{-(j-1)}\EE(r_0)=2\,\EE(r_0)\text{ for all } x\in K.
$$
Therefore $I(\mu)\lesssim\EE(r_0)<\infty$, and so $\Capacity(K)>0$.
\end{proof}
\end{lemma}

A more detailed and quantitative version of the previous proof may yield a lower bound for $\Capacity (K\cap B(x,r))$ for all $x\in K$ and $0<r<1$ which would imply that every point of $K$ is Wiener regular, even though $K$ does not satisfy the capacity density condition.

Kellogg's theorem \cite{Kellogg1928,Evans1933} (see also \cite[Theorem 6.32]{Prats-Tolsa-book}) asserts that the set of irregular points for the Laplace operator (hence for any operator $L_A$) has zero capacity. With this we are ready to prove \cref{lemma:positive elliptic measure}.

\begin{proof}[Proof of \cref{lemma:positive elliptic measure}]
By Kellogg's theorem and the fact that $\Capacity (K)>0$, we have that there exists a regular point $x \in K$. Let $\characteristic_K$ be the characteristic function of $K$, which is continuous on $\partial\Omega$. By Perron's method, $h_{\characteristic_K}^A=\omega^\cdot_A (K)$ is the $L_A$-harmonic extension in $\Omega$ of $\characteristic_K$. By the Wiener regularity of $x\in K$, $h_{\characteristic_K}^A(y) \to 1$ as $\Omega\ni y \to x$, whence we get that there is a point $y^\prime\in \Omega$ with $\omega^{y^\prime}_A(K)\geq 1/2$, say. Since $\Omega$ is connected, $\omega_A (K)>0$.
\end{proof}

\section{Proof of \texorpdfstring{\cref{thm:main}}{Theorem \getrefnumber{thm:main}}}\label{sec:proof}

In this section we prove \cref{thm:main}. We use the following standard fact relating the packing dimension of a set to the upper local dimension of a measure, see \cite[Proposition 2.3(c)]{Falconer1997} and also \cite[Chapter 6, Exercise 7]{Mattila1995}: if $\mu$ is a finite Borel measure on $\R^n$, $F\subset\R^n$ is a Borel set with $\mu(F)>0$, and $s\geq 0$ is such that
\begin{equation}\label{eq:upper local dimension}
\limsup_{r\to 0}\frac{\log\mu(B(x,r))}{\log r}\geq s \text{ for all } x\in F,
\end{equation}
then $\dim_P F\geq s$.

\begin{proof}[Proof of \cref{thm:main}]
    Let $K\subset B(0,1/4)$ and $\Omega=B(0,1)\setminus K$ be as in \cref{sec:construction}, and let $A$ be a uniformly elliptic matrix in $\Omega$ with ellipticity constant $\lambda$. Since the packing dimension of the elliptic measure does not depend on the pole, we fix $\omega_A\coloneqq\omega^p_{\Omega,A}$ for some $p\in\Omega$ with $|p|=1/2$.

    Let $E\subset\R^n$ with $\omega_A(E^c)=0$, and let $N$ be a Borel set with $E^c\subset N$ and $\omega_A(N)=0$. Then $F\coloneqq K\setminus N\subset E$ is a Borel set, and $\omega_A(F)=\omega_A(K)>0$ by \cref{lemma:positive elliptic measure}. Let $C\geq 1$ be the implicit constant in \cref{lemma:upper bound}, and let $x\in F$. For all $k$ large enough, since $0<\delta_k/2<1$,
    $$
    \frac{\log\omega_A(B(x,\delta_k/2))}{\log(\delta_k/2)} \geq \frac{\log C + \big(n-\frac{2}{k+1}\big)\log\delta_k}{\log\delta_k - \log 2} \xrightarrow{k\to\infty} n,
    $$
    because $\delta_k\to 0$. Hence \eqref{eq:upper local dimension} holds for $\mu=\omega_A$ and $s=n$, and so $\dim_P E\geq\dim_P F\geq n$. Taking the infimum over all such sets $E$, we get $\dim_P \omega_A = n$. Since $\Omega$ does not depend on $A$, this concludes the proof of the theorem.
\end{proof}

\subsection*{AI disclosure}

During the preparation of this paper, the author used AI to help with the verification of references, with checking computations, and with improving the wording of the text. The author reviewed and verified all the content and takes full responsibility for it.

\bibliographystyle{amsalpha}
\bibliography{references}

\providecommand{\bysame}{\leavevmode\hbox to3em{\hrulefill}\thinspace}
\providecommand{\MR}{\relax\ifhmode\unskip\space\fi MR }
\providecommand{\MRhref}[2]{%
  \href{http://www.ams.org/mathscinet-getitem?mr=#1}{#2}
}
\providecommand{\href}[2]{#2}
\begin{thebibliography}{GMPT25}

\bibitem[ABG26]{Albert-Badger-Genschaw-2026}
Michael Albert, Matthew Badger, and Alyssa Genschaw, \emph{{O}ne-ring
  {B}ourgain alternative and a 5-decimal bound on the dimension of harmonic
  measure}, \href{https://arxiv.org/abs/2610.28037}{arXiv:2610.28037}, 2026.

\bibitem[BG24]{Badger2024}
Matthew Badger and Alyssa Genschaw, \emph{Lower bounds on {B}ourgain's constant
  for harmonic measure}, Canad. J. Math. \textbf{76} (2024), no.~6, 1967--1986.
  \MR{4824556}

\bibitem[Bis92]{Bishop1992}
Christopher~J. Bishop, \emph{Some questions concerning harmonic measure},
  Partial differential equations with minimal smoothness and applications
  ({C}hicago, {IL}, 1990), IMA Vol. Math. Appl., vol.~42, Springer, New York,
  1992, pp.~89--97. \MR{1155854}

\bibitem[Bou87]{Bourgain1987}
J.~Bourgain, \emph{On the {H}ausdorff dimension of harmonic measure in higher
  dimension}, Invent. Math. \textbf{87} (1987), no.~3, 477--483. \MR{874032}

\bibitem[Car85]{Carleson1985}
Lennart Carleson, \emph{On the support of harmonic measure for sets of {C}antor
  type}, Ann. Acad. Sci. Fenn. Ser. A I Math. \textbf{10} (1985), 113--123.
  \MR{802473}

\bibitem[Cho04]{Choi2004}
Sunhi Choi, \emph{The lower density conjecture for harmonic measure}, J. Anal.
  Math. \textbf{93} (2004), 237--269. \MR{2110330}

\bibitem[CJT26a]{Cai-Jeznach-Tolsa-close-to-1}
Yingying Cai, Cole Jeznach, and Xavier Tolsa, \emph{{F}ailure of dimension drop
  for harmonic measure on {A}hlfors regular boundaries of dimension close to
  $1$}, \href{https://arxiv.org/abs/2610.15759}{arXiv:2610.15759}, 2026.

\bibitem[CJT26b]{Cai-Jeznach-Tolsa-flat}
\bysame, \emph{{F}ailure of dimension drop for harmonic measure on {A}hlfors
  regular flat boundaries},
  \href{https://arxiv.org/abs/2610.15805}{arXiv:\allowbreak 2610.15805}, 2026.

\bibitem[Eva33]{Evans1933}
G.~C. Evans, \emph{Application of {P}oincar\'{e}'s sweeping-out process}, Proc.
  Nat. Acad. Sci. U.S.A. \textbf{19} (1933), no.~4, 457--461.

\bibitem[Fal97]{Falconer1997}
Kenneth Falconer, \emph{Techniques in fractal geometry}, John Wiley \& Sons,
  Ltd., Chichester, 1997. \MR{1449135}

\bibitem[GM26]{Guillen-quasiconformal}
Ignasi Guill\'{e}n-Mola, \emph{The {H}ausdorff dimension of planar elliptic
  measures via quasiconformal mappings}, Rev. Mat. Iberoam. (2026), published
  online first.

\bibitem[GMPT25]{Guillen-Prats-Tolsa}
Ignasi Guill\'{e}n-Mola, Mart\'{\i} Prats, and Xavier Tolsa, \emph{The
  dimension of planar elliptic measures arising from {L}ipschitz matrices in
  {R}eifenberg flat domains}, Anal. Math. Phys. \textbf{15} (2025), no.~4,
  Paper No. 106, 79. \MR{4931423}

\bibitem[GW82]{Gruter1982}
Michael Gr\"{u}ter and Kjell-Ove Widman, \emph{The {G}reen function for
  uniformly elliptic equations}, Manuscripta Math. \textbf{37} (1982), no.~3,
  303--342. \MR{657523}

\bibitem[HK07]{Hofmann2007}
Steve Hofmann and Seick Kim, \emph{The {G}reen function estimates for strongly
  elliptic systems of second order}, Manuscripta Math. \textbf{124} (2007),
  no.~2, 139--172. \MR{2341783}

\bibitem[HKM06]{Heinonen2006}
Juha Heinonen, Tero Kilpel\"{a}inen, and Olli Martio, \emph{Nonlinear potential
  theory of degenerate elliptic equations}, Dover Publications, Inc., Mineola,
  NY, 2006, Unabridged republication of the 1993 original. \MR{2305115}

\bibitem[Hou95]{Housworth1995}
Elizabeth Housworth, \emph{Packing measure analysis of harmonic measure}, J.
  London Math. Soc. (2) \textbf{52} (1995), no.~3, 529--533. \MR{1363818}

\bibitem[JW88]{Jones1988}
Peter~W. Jones and Thomas~H. Wolff, \emph{Hausdorff dimension of harmonic
  measures in the plane}, Acta Math. \textbf{161} (1988), no.~1-2, 131--144.
  \MR{962097}

\bibitem[Kel28]{Kellogg1928}
O.~D. Kellogg, \emph{Unicit\'{e} des fonctions harmoniques}, C. R. Acad. Sci.
  Paris \textbf{187} (1928), 526--527.

\bibitem[KN85]{Kenig1985}
Carlos~E. Kenig and Wei-Ming Ni, \emph{On the elliptic equation {$Lu-k+K\,{\rm
  exp}[2u]=0$}}, Ann. Scuola Norm. Sup. Pisa Cl. Sci. (4) \textbf{12} (1985),
  no.~2, 191--224. \MR{829052}

\bibitem[LSW63]{Littman-Stampacchia-Weinberger-1963}
W.~Littman, G.~Stampacchia, and H.~F. Weinberger, \emph{Regular points for
  elliptic equations with discontinuous coefficients}, Ann. Scuola Norm. Sup.
  Pisa Cl. Sci. (3) \textbf{17} (1963), 43--77. \MR{161019}

\bibitem[LVV05]{Lewis-Verchota-Vogel-2005}
John~L. Lewis, Gregory~C. Verchota, and Andrew~L. Vogel, \emph{Wolff
  snowflakes}, Pacific J. Math. \textbf{218} (2005), no.~1, 139--166.
  \MR{2224593}

\bibitem[Mak85]{Makarov1985}
N.~G. Makarov, \emph{On the distortion of boundary sets under conformal
  mappings}, Proc. London Math. Soc. (3) \textbf{51} (1985), no.~2, 369--384.
  \MR{794117}

\bibitem[Mat95]{Mattila1995}
Pertti Mattila, \emph{Geometry of sets and measures in {E}uclidean spaces},
  Cambridge Studies in Advanced Mathematics, vol.~44, Cambridge University
  Press, Cambridge, 1995, Fractals and rectifiability. \MR{1333890}

\bibitem[{\O}ks81]{Oksendal1981}
Bernt {\O}ksendal, \emph{Brownian motion and sets of harmonic measure zero},
  Pacific J. Math. \textbf{95} (1981), no.~1, 179--192. \MR{631668}

\bibitem[{\O}ks88]{Oksendal1988}
\bysame, \emph{Dirichlet forms, quasiregular functions and {B}rownian motion},
  Invent. Math. \textbf{91} (1988), no.~2, 273--297. \MR{922802}

\bibitem[PT26]{Prats-Tolsa-book}
Mart{\'\i} Prats and Xavier Tolsa, \emph{{H}armonic measure in {E}uclidean
  spaces}, draft manuscript,
  \url{https://mat.uab.cat/~xtolsa/mesuraharmonica.pdf}, 2026, version of July
  28, 2026.

\bibitem[Swe92]{Sweezy1992}
Caroline Sweezy, \emph{The {H}ausdorff dimension of elliptic measure---a
  counterexample to the {O}ksendahl conjecture in {${\bf R}^2$}}, Proc. Amer.
  Math. Soc. \textbf{116} (1992), no.~2, 361--368. \MR{1161401}

\bibitem[Swe94]{Sweezy1994}
\bysame, \emph{The {H}ausdorff dimension of elliptic and elliptic-caloric
  measure in {${\bf R}^N,\;N\geq 3$}}, Proc. Amer. Math. Soc. \textbf{121}
  (1994), no.~3, 787--793. \MR{1186138}

\bibitem[Wol93]{Wolff1993}
Thomas~H. Wolff, \emph{Plane harmonic measures live on sets of
  {$\sigma$}-finite length}, Ark. Mat. \textbf{31} (1993), no.~1, 137--172.
  \MR{1230270}

\bibitem[Wol95]{Wolff1995}
\bysame, \emph{Counterexamples with harmonic gradients in {${\bf R}^3$}},
  Essays on {F}ourier analysis in honor of {E}lias {M}. {S}tein ({P}rinceton,
  {NJ}, 1991), Princeton Math. Ser., vol.~42, Princeton Univ. Press, Princeton,
  NJ, 1995, pp.~321--384. \MR{1315554}

\end{thebibliography}

\end{document}